\documentclass[a4paper]{amsart}
\usepackage[english]{babel}
\usepackage{amssymb}
\usepackage{graphicx}
\usepackage{booktabs}
\usepackage{longtable}
\usepackage[latin1]{inputenc}

\newcommand{\Z}{\mathbb{Z}}

\newcommand{\Q}{\mathbb{Q}}

\theoremstyle{definition}

\theoremstyle{remark}
\newtheorem*{rem}{Remark}
\begin{document}
\title[Elliptic curves and torsion]{Elliptic curves with torsion group $\Z/ 8 \Z $   or   $\Z/2 \Z \times \Z /6 \Z $}
\author[A. Dujella]{Andrej Dujella}
\address{Department of Mathematics\\University of Zagreb\\Bijeni{\v c}ka cesta 30, 10000 Zagreb, Croatia}
\email[A. Dujella]{duje@math.hr}
\author[J. C. Peral]{Juan Carlos Peral}
\address{Departamento de Matem\'aticas\\Universidad del Pa\'is Vasco\\Aptdo. 644, 48080 Bilbao, Spain}
\email[J. C. Peral]{juancarlos.peral@ehu.es}
\thanks{A. D. was supported by the Ministry of Science, Education and Sports, Republic of Croatia,
grant 037-0372781-2821. J. C. P. was supported by the UPV/EHU grant EHU 10/05.}
\subjclass[2010]{11G05}
\begin{abstract}
We show the existence of families of elliptic curves over $\Q$ whose generic rank is at least $2$ for the torsion groups $\Z/ 8 \Z $ and $\Z/2\Z \times \Z/6\Z$.  Also in  both cases
 we prove the existence of infinitely many elliptic curves, which are parameterized by the points of an elliptic curve with positive rank, with such torsion group and rank at least $3$.

These results  represent an improvement of previous results by Campbell, Kulesz, Lecacheux, Dujella and Rabarison where  families with rank at least $1$ were constructed in both cases.
\end{abstract}

\maketitle

\section{Introduction }

The construction of families of elliptic curves having high rank often is based in two basic strategies as mentioned by
 Elkies  in \cite{El}.

 a) The  Neron method studies the pencil of cubics passing through a set of nine rational random points and then looks for independence. See \cite{Sh} for a good description of the method. Families of rank up to $10$ where constructed in this way.

b) The Mestre  method  uses polynomial identities forcing the existence of rational points in the curve and then searchs for independence conditions.  In this way Mestre was able to construct a rank $11$ curve over $\Q(u)$, see \cite{Me}.

In our case we want that the  curve has a predetermined torsion group, $\Z/ 8 \Z $ or $\Z/2\Z \times \Z/6\Z$, so we should start from the general model for such torsion and then try to impose the existence of new points. One way to do this is looking for good quadratic sections. Also
the method developed by Lecacheux  \cite{Le1},  \cite{Le2}, can be used  in these cases. In her method  Lecacheux uses fibrations  of the corresponding surfaces, (as the ones given explicitlly  in Beauville \cite{Be}, in  Bertin and Lecacheux  \cite{BL} or in Livne and Yui  \cite{LY}) in order to find elliptic curves  with positive rank over $\Q(u)$. Another useful tool in order to get high rank curves over $\Q(u)$  are the diophantine triples, in fact for the torsion group $\Z/2\Z \times \Z/4\Z$  we construct  a rank $4$ family, the current record for this torsion group, see \cite{DP}.

In the case of torsion group $ \Z/8\Z$ families having    rank at least $1$ over $\Q(u)$ has been found by several authors, see  \cite{Ku}, \cite{Le1}, \cite{Le2} and \cite{Ra}.
In this paper we  prove the existence of two elliptic curves having  this torsion group
 and rank at least $2$ over $\Q(u)$,
also  the existence of infinitely many elliptic curves over $\Q$ with this torsion group  and rank at least $3$, parametrized by the points of an elliptic curve with positive rank.

Elliptic curves with torsion group $\Z/2\Z \times \Z/6\Z$ and rank at least $1$ over $\Q(u)$ have been constructed by several authors,  see  \cite{Ca}, \cite{Ku}, \cite{Le1}, \cite{Du1} and \cite{Ra}.
Here we prove the existence of three elliptic curves with  torsion group
$\Z/2\Z \times \Z/6\Z$ and whose rank over $\Q(u)$ is at least $2$.
We prove  also the existence of infinitely many elliptic curves with this torsion group   and rank at least $3$ over $\Q$, parametrized by the points of an elliptic curve of positive rank.

We describe first the model for the  elliptic curves having  each of these torsion groups. Then, in both cases, we show the existence of several  families having rank at least $1$ over $\Q(u)$. We finally show the existence of
two families in the $ \Z/8\Z$ case and three families in the case of torsion $\Z/2\Z \times \Z/6\Z$ having  rank at least $2$ over $\Q(u)$. In all the cases we describe the coefficients of the  families and the coordinates of  independent points.

We also exhibit some examples of curves with high rank.
The current records both for families and for curves can be found in \cite{Du2}.
One way to  find high rank elliptic curves over $\Q$ is the construction of elliptic curves with positive rank over $\Q(u)$, as high as possible, and then searching for good specializations  with adequate sieving  tools such as the Mestre-Nagao sums for example.

\section{Torsion group}

\subsection{Curves with  torsion group $\Z/8\Z $}
The Tate's normal form for an elliptic curve is given by
\[
E(b,c): \quad y^2+(1-c) x y - b y = x^3- b x^2
\]
(see \cite{Kn}). It is nonsingular if and only if $b\neq 0$.
  Using the adition law for $P=(0,0)$  and taking $d=b/c$  we have
\begin{align*}
4 P=&( d(d-1), d^2( c-d+1)), \\
- 4 P=&( d(d-1), d( d-1)^2)
\end{align*}
so $P$  is a torsion point of order $8$ for $b$ and $c$ as follows
\begin{align*}
b=& (2 v-1)(v-1),\\
c=&\frac{(2 v-1)(v-1)}{v}
\end{align*}
with $v$ a  rational, see \cite{Kn}.
For these values of $b$ and $c$ we can write the curve in the form $y^2= x^3+ A_{8}(v) x^2 + B_{8}(v) x$  where
\begin{align*}
A_{8}(v)=&1 - 8 v + 16 v^2 - 16 v^3 + 8 v^4, \\
B_{8}(v)=&16 (-1 + v)^4 v^4.
\end{align*}
Writing the curve in this form   is a convenient way  to search for candidates for new rational points. In fact  their $x$-coordinates  should be either divisors of $B$ or rational  squares times divisors of $B$.

\subsection{Curves with  torsion group $\Z/2\Z \times \Z/6\Z$}
Using again the Tate's normal form and the adition law  for  the point   $P=(0,0)$  we find
\begin{align*}
3 P&=(c, b-c), \\
- 3 P&=(c, c^2).
\end{align*}
It follows that $ P$ will be a torsion point of   order $6$  for  $b=c+c^2$.
For this value of $b$ we write the curve in the form $y^2= x^3+ A_6(c) x^2 + B_6(c) x$.  We get
\begin{align*}
A_6(c) &=-1 + 6 c - 3 c^2, \\
B_6(c)  &=-16 c^3.
\end{align*}
In the new coordinates the torsion point of order $6$ becomes $ (- 4 c, 4 c(1+c) )$.
 Now we use the fact that the curves with torsion $\Z/2\Z \times \Z/2\Z$ have a model $y^2=x(x-m)(x-n)$. So in order to get a curve with torsion group $\Z/2\Z \times \Z/6\Z$ it is enough that in the family above the polynomial
 \[
x^3+ A_6(c) x^2 + B_6(c) x=x(-16 c^3 + x + 6 c x - 3 c^2 x + x^2)
\]
 factorizes into linear  factors. So the discriminant $\Delta = (1 + c)^3 (1 + 9 c)$ of the second order polynomial must be a square, hence
 \[
c=\frac{-v^2+1}{2(3 v-5)}.
\]
For these values of $c$ the corresponding  curves have torsion group  $\Z/2\Z \times \Z/6\Z$ and can be written as
$y^2=x^3+ A_{26}(v) x^2+ B_{26}(v) x$ where
\begin{align*}
A_{26}(v) &=37 - 84 v + 102 v^2 - 36 v^3 - 3 v^4, \\
B_{26}(v)  &=32 (-1 + v)^3 (1 + v)^3 (-5 + 3 v).
\end{align*}

The torsion point of order $6$ transforms into
\[
(8 (-1+v) (1+v) (-5+3 v),8 (-3+v)^2 (-1+v) (1+v) (-5+3 v)).
\]

\begin{rem}
The cases treated here, jointly with the general curve having torsion $\Z/7\Z$  and whose model is
\[
y^2=x^3 + x^2 (1 - 2 t + 3 t^2 + 6 t^3 + t^4) +
 x (-8 t^2 (1 + t) (-1 + t + t^2)) + 16 t^4 (1 + t)^2,
 \]
 are the three cases in which the general model for such torsion group is a $K3$ surface.
\end{rem}

\section{Rank $1$ families}

\subsection{The case of torsion group $\Z/8\Z $}\label{z8r1}
For this torsion group we show ten conditions  upon $v$ leading to rank $1$ families. Some of them were already known due to the authors quoted before.
We first list eigth  values of $x_i$ which becomes the  $x$-coordinate of a new point once we specialize to the corresponding values of $v_i$, $i=1\ldots, 8$. We  include another two values, $v_9$ and $v_{10}$, found by Lecacheux  by using adequate fibration of the general model with torsion  $\Z/8\Z $.

\begin{align*}
x_1 &=\frac{- 16 v^4(1 - 4 v + 2 v^2)}{(-1 + 4 w)^2},&v_1=\frac{1+ w^2}{3 - 2 w + w^2},\\
x_2 &=\frac{-(-1 + v)^4 (-5 + 8 v) (-5 + 18 v)}{4 (-2 + 3 v)^2},&v_2=\frac{5 (1 + w^2)}{2 (9 + 4 w^2)},\\
x_3 &=\frac{-4 (-3 + v) (-1 + v)^2 v^4 (-1 + 3 v)}{(1 - 4 v + 2 v^2)^2},&v_3=\frac{1 + 3 w^2}{3 + w^2},\\
x_4 &=16 (-1 + v)^2 v^2 (1 - 2 v + 2 v^2),&v_4=\frac{(-2 + w) w}{-2 + w^2},\\
x_5 &=\frac{-64 (-1 + v)^2 v^2 (-1 - v + v^2)}{(-1 - 4 v + 4 v^2)^2},&v_5=\frac{(-2 + w) w}{1 + w^2},\\
x_6 &= -(-1 + v)^2 (1 - 6 v + 4 v^2),&v_6=\frac{2 - 2 w + w^2}{4+w^2},\\
x_7 &= 4 v^4,&v_7=\frac{-5+ w^2}{4(1+ w)},\\
x_8 &= \frac{-(-1 + v)^2 (-5 + 2 v)^2 (25 - 70 v + 36 v^2)}{(-7 + 6 v)^2},&v_8=\frac{34 - 6 w + w^2}{36 + w^2}.\\
\end{align*}
\[
v_9=\frac{w^2+ 12}{2(w^2+4)},\quad v_{10}=\frac{-2 w}{1-w+w^2}.
\]

Details of one of this cases is given  in the next section.
In every case the new point is of infinite order so the rank of the corresponding curve is at least $1$ over $\Q(w)$.

\subsection{The case of torsion group $\Z/2\Z \times \Z/6\Z$}\label{z26r1}

We have found several conditions  upon $v$ leading to rank $1$ families.
In the next table we list nine values of $x_i$ which becomes the  $x$-coordinate of a new point once we specialize to the  value of $v_i, i=1\ldots,9$.  As before some of them were already known due to the authors quoted above.

\begin{align*}
x_1 &=8 (-1+v)^3 (1+v),&v_1=\frac{3 (-1+w) (1+w)}{-29-8 w+w^2},\\
x_2 &=4 (1+v)^3,&v_2=\frac{3 (-3+w) (3+w)}{-45-24 w+w^2},\\
x_3 &=2 (-1+v) (1+v)^2 (-5+3 v),&v_3=\frac{-7+w^2}{1-4 w+w^2},\\
x_4 &=-16 (-1+v)^2 (1+v),&v_4=\frac{-11+w^2}{5-4 w+w^2},\\
x_5 &=16 (-5 + 3 v)(3v-7)^2,&v_5=\frac{3 (261+w^2)}{153-24 w+w^2},\\
x_6 &=16 (1 + v)(v-5)^2,&v_6=\frac{135-w^2}{141+24 w+w^2},
\end{align*}
\begin{align*}
x_7 &=\frac{4 (-1 + v)^2 (1 + v)^2 (41 - 54 v + 49 v^2)}{(-1 + 3 v)^2},&v_7=\frac{41-w^2}{2 (27+7 w)},\\
x_8 &=(-5 + 3 v) (3 v - 1)^2,&v_8=\frac{3}{5-w^2},\\
x_9 &=\frac{2(v-1)(v+1)^3(3v-1)^2}{(2 v+2)^2},&v_9=\frac{-7-2 w^2}{3(-3+ 2 w^2)}.\\
\end{align*}
More details  are given  in another section.
In every case the new point is of infinite order so the rank of the  curve is at least $1$ over $\Q(w)$.

\section{Rank $2$ families  for the torsion group  $\Z/8\Z$}

\subsection{A family with rank $1$}
We present here some details for the family  in which we have found two subfamilies with torsion  $\Z/8\Z$ and
generic rank at least $2$. It corresponds to the third  entry in the table of rank $1$ families above and it is a reparametrization of one of the families in \cite{Le1}.  By inserting in the general family $y^2= x^3+  A_{8}(v) x^2+  B_{8}(v) x$ the value $v=v_3(w)$ we get  the rank $1$ family given by $y^2= x^3+ AA_{8}(w) x^2+  BB_{8}(w) x$ where
\begin{align*}
AA_{8}(w)=&-31 - 148 w^2 + 214 w^4 - 116 w^6 + 337 w^8,\\
BB_{8}(w)=&256 (-1 + w)^4 (1 + w)^4 (1 + 3 w^2)^4.
\end{align*}

By searching on several  homogeneous spaces of the associate curve we have  found the possibility of imposing two new conditions which lead to new points. The  values $x_1$ and $x_2$ jointly with the specialization of the parameter are
\begin{align*}
x_1 &=\frac{(-1 + w)^2 (1 + w)^2 (5 + 7 w^2)^2 (11 + 25 w^2)}{16},&w_1=\frac{11-u^2}{10 u},\\
x_2 &=\frac{(-1 + w)^2 (1 + w)^2 (1 + 11 w^2)^2 (7 + 29 w^2)}{16 w^2},&w_2=\frac{29 - 12 u + u^2}{-29 + u^2}.
\end{align*}

 With these specializations we get  two different families of rank at least $2$ over $\Q(u)$.

 \subsection{First family with rank $2$ and torsion group $\Z/8\Z$}
 Once we insert $w_1$ into the coefficients $ AA_{8},BB_{8}$ we get as new coefficients  $AAA_{8}, BBB_{8}$ given by
\begin{align*}
AAA_{8} &=   337 u^{16} - 41256 u^{14}  + 4047356 u^{12}  -288332632 u^{10 } + 2363813190 u^8 \\
& - 34888248472 u^6 + 59257339196 u^4- 73087520616 u^2+72238942897, \\
BBB_{8} &= 256\, (363 + 34 u^2 + 3 u^4)^4\,(11 + u)^4\,(-11 + u)^4\,(-1 + u)^4\, (1 + u)^4.
\end{align*}
The $x$-coordinates of  two  independent infinite  order points   are

\begin{align*}
X_{1}&=\\
 &\frac{2^{12} 5^2 (-11 + u)^2 (-1 + u)^2 u^2 (1 + u)^2 (11 + u)^2 (-11 +
   u^2)^2 (363 + 34 u^2 + 3 u^4)^4}{(102487 - 303468 u^2 + 43482 u^4 - 2508 u^6 + 7 u^8)^2},  \\
X_{2} &= \\
&\frac{(-11 + u)^2 (-1 + u)^2 (1 + u)^2 (11 + u)^2 (11 + u^2)^2 (847 +
   346 u^2 + 7 u^4)^2}{64 u^2}.
\end{align*}

The $x$-coordinates of  the torsion point of order $8$ is:

\begin{align*}
T_{1}=
 &-8 (-11 + u) (-1 + u) (1 + u) (11 + u) (363 + 34 u^2 + 3 u^4)^3.
\end{align*}
That the rank of this curve is at least $2$ over $\Q(u)$ can be proved using a specialization argument,
since the specialization map is a homomorphism.

 \subsection{Second family with rank $2$ and torsion group $\Z/8\Z$}

Once we insert $w_2$ into the coefficients $ AA_{8},BB_{8}$ we get as new coefficients  $aaa_{8}, bbb_{8}$ given by
\begin{align*}
aaa_{8} &=   500246412961 - 2069985157080 u + 3162080774436 u^2 -
 2895517882032 u^3 +\\& 1873181389706 u^4 - 906769167048 u^5 +
 333391978480 u^6 - 93284915496 u^7 +\\& 19860033555 u^8 -
 3216721224 u^9 + 396423280 u^{10} - 37179432 u^{11} +\\& 2648426 u^{12} -
 141168 u^{13} + 5316 u^{14} - 120 u^{15} + u^{16}, \\
bbb_{8} &=256 (-6 + u)^4 u^4 (-29 + 6 u)^4 (841 - 522 u + 137 u^2 - 18 u^3 +
   u^4)^4.
\end{align*}

The $x$-coordinates of  two  independent infinite  order points   are

\begin{align*}
&\\
 &\frac{64 (-6 + u)^2 u^2 (-29 + 6 u)^2 (-29 + u^2)^2 (29 - 12 u +
   u^2)^2 (841 - 522 u + 137 u^2 - 18 u^3 + u^4)^4}{(707281 - 292668 u - 200158 u^2 + 168432 u^3 - 46685 u^4 + 5808 u^5 -
  238 u^6 - 12 u^7 + u^8)^2},  \\
& \\
&\frac{(-6 + u)^2 u^2 (-29 + 6 u)^2 (87 - 29 u + 3 u^2)^2 (2523 - 1914 u +
   541 u^2 - 66 u^3 + 3 u^4)^2}{4 (29 - 12 u + u^2)^2}.
\end{align*}
The  $x$-coordinate of  the torsion point of order $8$ is:
\begin{align*}
T_{2}=
 &8 (-6 + u) u (-29 + 6 u) (841 - 522 u + 137 u^2 - 18 u^3 + u^4)^3.
\end{align*}

\begin{rem}
When we use $w_3=\frac{(-3 + u) (3 + u)}{7-6 u}$ and $w_4=\frac{(-3 + u) (3 + u)}{11-6 u}$ in the family of rank  at least $1$ corresponding to $v_{10}$ we get two families of rank at least $2$.  They are a reparametrization of the two families above.
So we have the following fact, when we specialize in the general family  with torsion group $\Z/8\Z$ to
$v_3=\frac{1+ 3 w^2}{3+ w^2}$ and to $v_{10}=\frac{-2w}{1-w+w^2}$ we get two diferent  families having rank at least $1$. When we use in the first the values $w_1=\frac{11- u^2}{10 u}$ and $w_2=\frac{29-12 u+ u^2}{-29+ u^2}$  we get two families of rank at least $2$ that are a reparametrization of the rank $2$ families that we get by using $v_{10}$ followed by  $w_3=\frac{(-3+u)(3+u)}{7-6u}$ and $w_4=\frac{(-3+u)(3+u)}{11-6u}$.

So at the end with the changes $v_3$ follow by $w_1$ and $w_2$ and $v_{10}$ follow by $w_3$ and $w_4$ we reach the same families of rank at least $2$ even though the families of rank at least $1$ of the intermediate step are different.
\end{rem}

\subsection{Rank $3$ for the  torsion group $\Z/8\Z$}
It can be proved that there exist infinitely many elliptic curves with torsion group $\Z/ 8 \Z $ parametrized by the points of a positive rank elliptic curve.
In fact it  is enough to see that the equation $w_1(r)=w_2(s)$, i.e.:
\[
\frac{11-r^2}{10 r}=\frac{29 - 12 s + s^2}{-29 + s^2}
\]
  has infinitely many solutions. This is the same as to solve
  \[
  319 + 290 r - 29 r^2 - 120 r s - 11 s^2 + 10 r s^2 + r^2 s^2=0
  \]
   in rational terms, so   the discriminant  $\Delta =3509 + 62 r^2 + 29 r^4$ has to be a square.
   But  $t^2=3509 + 62 r^2 + 29 r^4$ has a solution, $(r,t)=(1,60)$ for example,
   hence it is equivalent to the cubic   $y^2= x^3 - 463 x^2+ 45936 x$
   whose rank is $2$ as proved with {\tt mwrank} \cite{Cr}.
   This, jointly with the independence of the corresponding points, implies the  existence of infinitely many solutions parametrized by the points of the elliptic curve,
   see \cite{Le1} or \cite{Ra} for this kind or argument.

\section{Rank $2$ families  for the torsion group $\Z/2\Z \times \Z/6\Z$}

\subsection{A family with rank $1$}
We present here some details for the family of rank $1$  in which we have found two subfamilies of
generic rank at least $2$. It corresponds to the eighth entry in the table above.
We have found the possibility of converting $x_8=(-5 + 3 v) (3 v - 1)^2$ in the $x-$coordinate of a new point by considering the homogenous space $(U,V)=(3 v-1, 1)$ of the initial family
\[
y^2=x^3+ (37 - 84 v + 102 v^2 - 36 v^3 - 3 v^4) x^2+ (32 (-1 + v)^3 (1 + v)^3 (-5 + 3 v)) x.
\]
The condition that has to be  fulfilled is that $  v(-3+ 5v)$ converts into a square, hence  we get  $v_8=\frac{3}{5-w^2}$.
Once we insert $v_8$ in the preceding family and take  off denominators we have the   family  $y^2=x^3+ AA_{26}(w) x^2+ BB_{26}(w) x$ where
\begin{align*}
AA_{26}(w) &=9472 - 7808 w^2 + 2688 w^4 - 488 w^6 + 37 w^8, \\
BB_{26}(w)  &=32 (-8 + w^2)^3 (-5 + w^2) (-2 + w^2)^3 (-16 + 5 w^2).
\end{align*}
The point of infinite order is
\begin{align*}
P=&\Big{ (}-(-5+w^2) (4+w^2)^2 (-16+5 w^2),\\
&27  (-2+w)^2 w (2+w)^2 (-5+w^2)
(4+w^2) (-16+5 w^2)\Big{)}
\end{align*}
and the torsion point of order 6 is
\begin{align*}
T=&\Big{(} 8 (-8+w^2) (-5+w^2) (-2+w^2) (-16+5 w^2),\\
&72 (-2+w)^2 (2+w)^2 (-8+w^2) (-5+w^2) (-2+w^2) (-16+5 w^2)\Big{)}
\end{align*}

By searching on several  homogeneous spaces we have  found the possibility of imposing two new conditions which lead to new points, hence, in case of independence, to a couple of rank $2$ families. The  values $x_1$ and $x_2$ jointly with the specialization of the parameter are

\begin{align*}
x_1 &=2 (-8+ w^2)^2(- 2+ w^2)^3,&w_1=\frac{2 (7 + u^2)}{-7 - 2 u + u^2},\\
x_2 &=(-8 + w^2) (- 5 + w^2) (-16 + 5 w^2) w^4,&w_2=\frac{5 - 2 u + u^2}{-5 + u^2}.
\end{align*}

\subsection{First  family with rank $2$}\label{z26r2} In order to force $2 (-8 + w^2)^2 (-2 + w^2)^3$ to be the $x$-coordinate of a new point  it is enough to solve   $2 w^2-7= M^2$. This is achieved with
$w_1=\frac{2(7+u^2)}{-7- 2 u+ u^2}$, and  the corresponding family is $y^2=x^3+ AAA_{26}(u) x^2+ BBB_{26}(u) x$ where
\begin{align*}
AAA_{26}(u)= &-2 (5764801 + 6588344 u - 21647416 u^2 + 29445864 u^3 - 9604 u^4 + \\
  & 27969592 u^5 - 44631944 u^6 + 9779112 u^7 + 5909830 u^8 -
   1397016 u^9\\& - 910856 u^{10} - 81544 u^{11} - 4 u^{12} - 1752 u^{13} -
   184 u^{14} - 8 u^{15} + u^{16}),\\
BBB_{26}(u) = &(-7 - 10 u + u^2)^3 (-7 + 2 u + u^2)^3 (49 + 140 u - 106 u^2 -
   20 u^3 + u^4)\\& (49 - 28 u + 38 u^2 + 4 u^3 + u^4)^3 (49 - 112 u +
   110 u^2 + 16 u^3 + u^4).
\end{align*}
The $x$-coordinates of the two infinite order points are
\begin{align*}
X_1=&(49 + 140 u - 106 u^2 - 20 u^3 + u^4) (49 + 14 u + 2 u^2 - 2 u^3 +
   u^4)^2\times\\& (49 - 112 u + 110 u^2 + 16 u^3 + u^4),\\
X_2=&\frac{(-7 - 10 u + u^2)^2  (-7 + 2 u + u^2)^2 (49 - 28 u + 38 u^2 + 4 u^3 + u^4)^3}{(-7 - 2 u + u^2)^2},
\end{align*}
and the $x$-coordinate of the torsion point of order $6$ is
\begin{align*}
T=&(-7 - 10 u + u^2) (-7 + 2 u + u^2) (49 + 140 u - 106 u^2 - 20 u^3 +
   u^4) \times \\& (49 - 28 u + 38 u^2 + 4 u^3 + u^4) (49 - 112 u + 110 u^2 + 16 u^3 +
   u^4).
\end{align*}

The specialization for $u=2$ gives an elliptic curve with rank $2$ and the specialized points are independent so the fact that specialization is a homomorphism  implies that the rank of the curve  $y^2=x^3+ AAA_{26}(u) x^2+ BBB_{26}(u) x$  is at least $2$ over $\Q(u)$.

\subsection{Second  family with rank $2$} Imposing $(-8 + w^2) (- 5 + w^2) (-16 + 5 w^2) w^4$ as the $x$-coordinate of a new point  it is equivalent to solve   $5 w^2-4= M^2$. This is achieved with
$w_2=\frac{5 - 2 u + u^2}{-5 + u^2}$.

The corresponding family is $y^2=x^3+ aaa_{26}(u) x^2+ bbb_{26}(u) x$ where
\begin{align*}
aaa_{26}(u)= &1523828125 + 1171250000 u - 3482125000 u^2 - 1970850000 u^3 + \\
  &
 3530367500 u^4 + 1221154000 u^5 - 2018502200 u^6 - 238418640 u^7 + \\
  &
 632792782 u^8 - 47683728 u^9 - 80740088 u^{10} + 9769232 u^{11} +\\&
 5648588 u^{12} - 630672 u^{13} - 222856 u^{14} + 14992 u^{15} + 3901 u^{16},\\
bbb_{26}(u) = &128 (-7 + 2 u + u^2)^3 (-25 - 10 u + 7 u^2)^3 (25 + 5 u - 16 u^2 +
   u^3 + u^4)\times \\& (25 + 20 u - 34 u^2 + 4 u^3 + u^4)^3 (275 + 100 u -
   230 u^2 + 20 u^3 + 11 u^4)
\end{align*}

The $x$-coordinates of the two infinite order points are
\begin{align*}
X_1=&-4 (25 + 5 u - 16 u^2 + u^3 + u^4) (125 - 20 u - 26 u^2 - 4 u^3 +
   5 u^4)^2\times\\& (275 + 100 u - 230 u^2 + 20 u^3 + 11 u^4),\\
X_2=&\frac{-4 (5 - 2 u + u^2)^4 (-7 + 2 u + u^2)(-25 - 10 u +
    7 u^2) }{(-5 + u^2)^2}\times\\& (25 + 5 u - 16 u^2 + u^3 + u^4) (275 + 100 u - 230 u^2 +
    20 u^3 + 11 u^4),
\end{align*}
and the $x$-coordinate of the torsion point of order $6$ is
\begin{align*}
T=&32 (-7 + 2 u + u^2) (-25 - 10 u + 7 u^2) (25 + 5 u - 16 u^2 + u^3 +
   u^4)\times\\& (25 + 20 u - 34 u^2 + 4 u^3 + u^4) (275 + 100 u - 230 u^2 +
   20 u^3 + 11 u^4).
\end{align*}

\subsection{Third  family with rank $2$}

A variant of the model for torsion group  $\Z/6\Z$   by Hadano (see \cite{H}) will be used  for the construction of another  family of curves with torsion group $\Z/2\Z \times \Z/6\Z$ and rank at least $2$.
The curves with torsion group $\Z/6\Z$ in \cite{H} have the equation
\[
Y^2=X^3+X^2(a^2+ 2 a b- 2 b^2)- X(2 a-b) b^3.
\]
For this model we first force the existence of a new point in order to have a family with this torsion and rank at least $1$, then we choose the parameters that give complete factorization,  and so the torsion $\Z/2\Z \times \Z/6\Z$, and in the resulting family for this torsion and rank $1$ we produce a second independent point with a quadratic section.

We observe that with $b=-\frac{(-1 + a^2 - v) (-1 + a^2 + v)}{4 (-a + a^3 - v)}$  we have a point in the curve with  $x$ coordinate  given by
\[
X=-\frac{(1 - 4 a + 3 a^2 - v) (-1 + a^2 + v)^3}{16 (-a + a^3 - v)^2}
\]
Once we clear denominators this family can be written as  $Y^2= X^3+ A_{61}(a,v) X^2+ B_{61}(a,v) X$ where

\begin{align*}
A_{61}&=2 (-1 + 8 a^2 - 10 a^4 + 3 a^8 + 4 a v + 8 a^3 v - 12 a^5 v + 2 v^2 +
   6 a^4 v^2 - 4 a v^3 - v^4), \\
B_{61}&=(-1 + a^2 - v)^3 (1 - 4 a + 3 a^2 - v) (-1 + a^2 + v)^3 (1 + 4 a +
   3 a^2 + v).
\end{align*}
The new point becomes
\[
P=(-(1 - 4 a + 3 a^2 - v) (-1 + a^2 + v)^3, 4 (1 - 4 a + 3 a^2 - v) (-a + a^3 - v) (-1 + a^2 + v)^3).
\]
For $a\neq 1$, in general, $P$ is of infinite order so the curve has rank at least $1$ over $\Q(a,v)$. Now we observe that the complete factorization of the cubic is equivalent to forcing the discriminant of a  second degree polynomial to be a square.

The condition is $a (-a + a^3 - v) (-1 + a^2 - a v + v^2)=$ square. This can be achieved with $a=v+1$ followed by  $v=\frac{1- w^2}{-3+ 2 w}$. Once we perform these changes we get  $Y^2=X^3+a_{26}(w) X^2+ b_{26}(w) X$ where
\begin{align*}
a_{26}&=2 (-24 - 216 w + 1008 w^2 - 1596 w^3 + 1319 w^4 - 648 w^5 + 198 w^6 -
   36 w^7 + 3 w^8), \\
b_{26}&=(-4 + w)^3 (-3 + w) (-2 + w)^3 (-1 + w)^3 w (1 + w)^3 (-7 +
   3 w) (-2 + 3 w).
\end{align*}
With these changes the $x$-coordinate of the infinite order point is
\[
X=-\frac{(-4 + w)^3 (-2 + w)^3 (-1 + w)^2 w (1 + w)^2 (-2 + 3 w)}{(2 - 2 w + w^2)^2}
\]
So we have a curve with torsion group $\Z/2\Z \times \Z/6\Z$ and rank at least $1$ over $\Q(w)$.

Now we see that $(-4 + w) (-3 + w) (-1 + w)^3 (1 + w)^2 (-7 + 3 w)$ will be a new point on the curve if we force  $4 - 9 w + 3 w^2$ to be a square. This is the same as choosing  $w=-\frac{9 + 4 u}{-3 + u^2}$.  Once we perform this change  and clear denominators we get the following coefficients for the new family
\begin{align*}
A_{263}=&-2 (157464 - 1889568 u - 13594392 u^2 - 38047968 u^3 -\\
& 62500248 u^4 -
   69622416 u^5 - 57719412 u^6 - 38941344 u^7 -\\
   & 23353995 u^8 -
   12980448 u^9 -
    6413268 u^{10} - 2578608 u^{11} -\\
    & 771608 u^{12} -
   156576 u^{13} - 18648 u^{14} - 864 u^{15} + 24 u^{16}), \\
B_{263}=&-(-6 + u)^3 u (2 + u)^3 (-1 + 2 u)^3 (3 + 2 u)^3 (4 + 3 u) (9 +
   4 u)\\
   & (6 + 4 u + u^2)^3 (3 + 4 u + 2 u^2)^3 (21 + 12 u +
   2 u^2) (6 + 12 u + 7 u^2).
\end{align*}
The $x$-coordinates of the non-torsion points are
\begin{align*}
X_1=&(-6 + u)^2 u (2 + u)^2 (-1 + 2 u) (3 + 2 u) (4 + 3 u)\times\\
& (6 + 4 u +
   u^2)^3 (6 + 12 u + 7 u^2),\\
   X_2=&-\frac{(-6 + u)^2 (2 + u)^2 (-1 + 2 u)^3 (3 + 2 u)^3 (9 + 4 u) }{(45 + 48 u + 22 u^2 + 8 u^3 + 2 u^4)^2}\times\\
   &(6 + 4 u +
   u^2)^2 (3 + 4 u + 2 u^2)^3 (21 + 12 u + 2 u^2).
\end{align*}
A specialization argument, as in Section \ref{z26r2}, shows that these two  points are independent so this  curve has rank at least $2$ over $\Q(u)$ and  torsion group $\Z/2\Z \times \Z/6\Z$.

\subsection{Rank $3$ for the  torsion group $\Z/2\Z \times \Z/6\Z$}
Now we prove the existence of infinitely many elliptic curves with rank at least $3$ and torsion group $\Z/2\Z \times \Z/6\Z$ parametrized by the points of a positive rank elliptic curve.
In fact it  is enough to see that the equation $w_1(r)=w_2(s)$, i.e.:
\[
\frac{2 (7 + r^2)}{-7 - 2 r + r^2}=\frac{5 - 2 s + s^2}{-5 + s^2}
\]
  has infinitely many solutions. This is the same as solving
  \[
 -35 + 10 r - 15 r^2 - 14 s - 4 r s + 2 r^2 s + 21 s^2 + 2 r s^2 +
 r^2 s^2=0
  \]
   in rational terms, so   the discriminant  $\Delta =49 - 7 r + 20 r^2 + r^3 + r^4$ has to be a rational square. But  $t^2=49 - 7 r + 20 r^2 + r^3 + r^4$ has a solution, $(r,t)=(1, 8)$ for example, hence it is birationally equivalent to the cubic   $y^2= x^3 - 43 x^2+ 280 x$ whose rank is $1$ as proved with {\tt mwrank} \cite{Cr}. This and the independence of the corresponding points, implies the  existence of infinitely many solutions parametrized by the points of the elliptic curve.

  \section{Examples of curves with high rank}
  \subsection{The case of torsion group $\Z/8\Z $}
         The highest known rank of an elliptic curve over $\Q$ with  torsion group $\Z/8\Z$ is rank $6$ curve found by Elkies in $2006$. See \cite{Du2} for the details of this curve.

The following list includes examples of rank $5$ curves  found in the rank $1$ families of subsection \ref{z8r1}. First column indicates the number of the  family and the second the value(s) of the parameter producing a rank $5$ curve. The indication (L)  means that this curve has been previously  found by Lecacheux.
\begin{align*}
 \hbox{Family number}&  &\hbox{$w$ values}\\
2 \hspace{1cm}&& \frac{287}{109},\\
3 \hspace{1cm}&& \frac{73}{83},\,\, \frac{37}{157},\\
4 \hspace{1cm}&&- \frac{87}{28},\\
5 \hspace{1cm}&& \frac{317}{10},\\
6 \hspace{1cm}&& -\frac{28}{79},\,\, \frac{100}{29} \,(\hbox{L}),\,\, \frac{304}{55}.\\
\end{align*}

       \subsection{The case of torsion group $\Z/2\Z \times \Z/6\Z$}

       The highest known rank of an elliptic curve over $\Q$ with  torsion group $\Z/2\Z \times \Z/6\Z$ was found by Elkies in $2006$. See \cite{Du2} for the details of this curve. It has rank $6$  and it can be seen to correspond to the value $u=-\frac{16}{3}$ in the family with rank $2$ included in subsection \ref{z26r2}. For $u=\frac{5}{13}$ we get a curve with  rank $5$.

 The following list includes examples of rank $5$ curves  found in the rank $1$ families of subsection \ref{z26r1}. First column indicates the number of the  family and the second the value(s) of the parameter producing a rank $5$ curve. The indication (L) and (D) means that in these cases the curves were found previously by Lecacheux and Dujella respectively.

       \begin{align*}
\hbox{Family number}&  &\hbox{$w$ values}\\
1 \hspace{1cm}&  &\frac{306}{11}\\
3 \hspace{1cm}&& \frac{53}{90} \,(\hbox{D}),\,\, -\frac{127}{74} \,(\hbox{L})\\
5 \hspace{1cm}&& \frac{31}{42}\\
6 \hspace{1cm}&& \frac{13}{43}, \,\, -\frac{431}{33}\\
7 \hspace{1cm}&& \frac{115}{6}, \,\, -\frac{142}{33}\, (\hbox{D}), \,\, -\frac{391}{387}, \,\, -\frac{1011}{551}\\
8 \hspace{1cm}&& \frac{302}{161}\\
9 \hspace{1cm}&& \frac{44}{61}, \,\, \frac{40}{57}, \,\, \frac{172}{191}, \,\, \frac{214}{163}, \,\, \frac{284}{197}
\end{align*}

 \end{document}